\documentclass[reqno]{amsart}

\usepackage[T1]{fontenc}
\usepackage{amsmath,amssymb,amsthm}
\usepackage{tikz}
\usepackage[colorlinks,linkcolor=black,citecolor=black,urlcolor=blue]{hyperref}

\theoremstyle{plain}
\newtheorem{theorem}{Theorem}[section]
\newtheorem{lemma}[theorem]{Lemma}
\newtheorem{proposition}[theorem]{Proposition}
\newtheorem{corollary}[theorem]{Corollary}

\newtheorem{mainthm}{Theorem}

\theoremstyle{definition}
\newtheorem{definition}[theorem]{Definition}
\newtheorem{example}[theorem]{Example}

\theoremstyle{remark}
\newtheorem{remark}[theorem]{Remark}

\numberwithin{equation}{section}

\newcommand{\cl}{\operatorname{cl}}
\newcommand{\col}{\operatorname{col}}
\newcommand{\wid}{\operatorname{width}}
\newcommand{\Z}{\operatorname{Z}}
\newcommand{\M}{\operatorname{M}}
\newcommand{\GP}[2]{P(#1,#2)}                    
\newcommand{\Pinf}[1]{P(\infty,#1)}              
\newcommand{\ZZ}{\mathbb{Z}}                     

\begin{document}

\title[Zero forcing in periodic graphs]
      {Stabilization of zero forcing numbers\\ of finite quotients of periodic graphs}

\author{Young H. Kim}
\address{Tuskegee University, Tuskegee, Alabama, USA}
\email{ykim@tuskegee.edu}

\author{Youngsoo Kim}
\address{Tuskegee University, Tuskegee, Alabama, USA}
\email{ykim1@tuskegee.edu}

\author{Byunghoon Lee}
\address{Tuskegee University, Tuskegee, Alabama, USA}
\email{blee@tuskegee.edu}


\begin{abstract}
  We give a method to determine the zero forcing number of a graph 
  by passing to an infinite cover of the graph. 
  When the quotient is a covering map, forcing commutes with taking preimages.
  If, in addition, the closure of every small set in the cover touches none of its translates,
  then no small set is a zero forcing set of the quotient.
  For periodic covers, this condition reduces to a property of the cover itself,
  and it gives an explicit threshold beyond which the zero forcing number of the quotients stabilizes 
  to that of the cover. 
  We apply the method to several classes of graphs: the generalized Petersen graphs,
  the Cartesian products with a cycle, and the circulant graphs.
\end{abstract}

\maketitle

\section{Introduction}
Let $G$ be a simple undirected graph with vertex set $V(G)$ and edge set $E(G)$.
An edge is written as $xy$ for vertices $x$ and $y$ and we write $x\sim y$ if $xy \in E(G)$.
We denote the neighbors of a vertex $x$ by $N(x)$, and assume throughout the paper
$|N(x)|<\infty$ for all vertices $x$. 
For any $X \subseteq V(G)$, the \emph{boundary} of $X$ is 
\[
\partial X := \{x \in X : \text{ there exists } y \in V(G) \setminus X \text{ such that } x \sim y \}.
\]

\begin{figure}[ht]
  \centering
  \begin{tikzpicture}[scale=0.86,
      bvertex/.style = {circle, fill=black, draw=black, line width=.5pt,
                        inner sep=0pt, minimum size=1.6mm},
      nvertex/.style = {circle, fill=black!35, draw=black, line width=.5pt,
                        inner sep=0pt, minimum size=1.6mm},
      wvertex/.style = {circle, fill=white, draw=black, line width=.5pt,
                        inner sep=0pt, minimum size=1.6mm},
      edge/.style    = {line width=.4pt, black!40},
      force/.style   = {line width=1pt, black, -stealth, shorten >=1mm}]
    \def\rounds{0/0/2, 1/0/1, 2/0/1, 3/0/1, 4/0/1, 5/0/1, 6/0/1,
                7/0/2, 8/3/2, 9/3/2, 10/4/3, 11/3/2, 12/3/2}
    \def\forces{1/u1/v1, 1/u2/v2, 1/u3/v3, 1/u4/v4, 1/u5/v5, 1/u6/v6,
                2/v3/v0, 2/v4/v7, 2/v5/v8, 2/v6/v9, 2/v1/v11, 2/v2/v12,
                3/u7/u8, 3/v9/u9, 3/v7/v10, 3/v11/u11, 3/u0/u12,
                4/u11/u10}
    \def\pickstyle#1{%
      \ifnum#1>\s \def\st{wvertex}%
      \else\ifnum#1=\s \ifnum\s>0 \def\st{nvertex}\else\def\st{bvertex}\fi
      \else \def\st{bvertex}\fi\fi}
    \foreach \s/\px/\py in {0/0/0, 1/4.8/0, 2/9.6/0, 3/2.4/-5.0, 4/7.2/-5.0} {
      \begin{scope}[shift={(\px,\py)}]
        \foreach \i in {0,...,12} {
          \coordinate (u\i) at ({90-360*\i/13}:1.65cm);
          \coordinate (v\i) at ({90-360*\i/13}:1.0cm);
        }
        \foreach \i in {0,...,12} {
          \pgfmathtruncatemacro{\j}{mod(\i+1,13)}
          \pgfmathtruncatemacro{\l}{mod(\i+3,13)}
          \draw[edge] (u\i) -- (u\j);
          \draw[edge] (u\i) -- (v\i);
          \draw[edge] (v\i) -- (v\l);
        }
        \foreach \r/\a/\b in \forces {
          \ifnum\r=\s \draw[force] (\a) -- (\b); \fi
        }
        \foreach \i/\su/\sv in \rounds {
          \pickstyle{\su}\node[\st] at (u\i) {};
          \pickstyle{\sv}\node[\st] at (v\i) {};
          \node[font=\tiny, text=black!55] at ({90-360*\i/13}:1.98cm) {\i};
        }
        \node at (0,-2.4) {$B_{\s}$};
      \end{scope}
    }
  \end{tikzpicture}
  \caption{The forcing process in $\GP{13}{3}$ from the set $S=B_0$ of 
  eight consecutive outer vertices.}
  \label{fig:P13-3}
\end{figure}
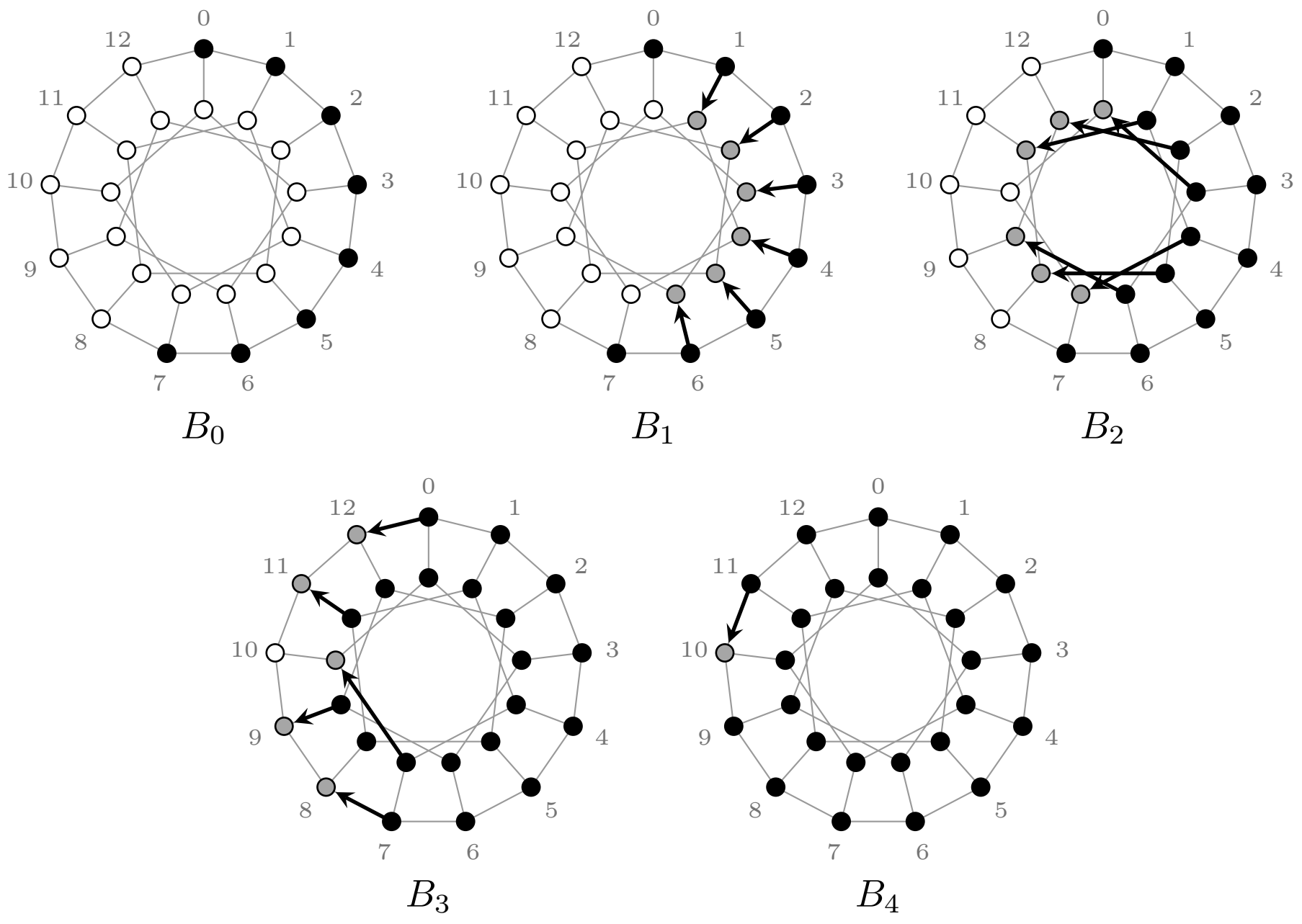

We play a game of coloring. Each vertex is colored either black or white.
Let $S$ be the initial set of black vertices of $G$. 
The color change rule says that if a black vertex $u$ has exactly one white neighbor $v$, 
then $v$ becomes black. 
We say $u$ \emph{forces} $v$, and call $u$ a \emph{forcer} of $v$.
A vertex $v$ may have several forcers, but each vertex can force at most once. 
Define recursively $B_0=S$ and $B_{i+1}=B_i\cup W_i$ where 
$W_i$ is the set of all unique white neighbors of some vertex in $B_i$.
\[
W_i = \{v \in V(G)\setminus B_i : \exists x\in B_i \text{ such that } N(x)\setminus B_{i}=\{v\} \}
\]
Then we define the \emph{closure} of $S$ by
\[
\cl(S) = \bigcup_{i=0}^\infty B_i.
\]
(See Figure \ref{fig:P13-3} for an example of the forcing process.)
A set $S$ is a \emph{zero forcing set} if $\cl(S) = V(G)$.
The smallest size of a zero forcing set of $G$ is called the \emph{zero forcing number}
of $G$ and is denoted by $\Z(G)$.

\begin{remark}
  By definition, $B_0\subseteq B_1\subseteq B_2\subseteq \cdots$. 
  For a finite graph the sets $B_i$ stop increasing at some point and $\cl(S) = B_n$ for some $n$.
  In the infinite graphs infinitely many rounds may be needed.
\end{remark}

Zero forcing was introduced in~\cite{AIM2008} as a combinatorial upper bound
for the maximum nullity $\M(G)$ over real symmetric matrices whose off-diagonal
nonzero pattern is given by $G$. One always has $\M(G)\le\Z(G)$.
The same process was introduced independently to study the controllability
of quantum systems~\cite{BurgarthGiovannetti2007}.
It has since become a graph parameter studied in its own right.
See~\cite{Barioli2010} for forcing chains and related parameters,
and~\cite{KalinowskiKamcevSudakov} for bounds on general classes of graphs.
The zero forcing number has been determined for many particular families of graphs.
These include paths, cycles, complete graphs, grids and other Cartesian products~\cite{AIM2008,Benson2018},
circulants~\cite{Meyer2012,DeaettMeyer2016,Duong2020},
gear and helm graphs~\cite{Anderton2023}, peony and web graphs~\cite{Anderton2024},
twisted hypercubes~\cite{CollierJanssen2025},
and generalized Petersen graphs~\cite{Alameda2018,RPT2020,Krishnan2026,KMZ2026}.
Upper bounds are determined by finding particular zero forcing sets.
For the lower bounds, these papers mostly give an argument tailored to the family,
construct a matrix of the right nullity (which also proves $\M(G)=\Z(G)$),
or check each graph by computer.

In this paper, we introduce a new technique to prove the lower bound on $\Z(G)$ for 
all but finitely many quotients of a single infinite graph at once.
Let $\widetilde{G}$ be a connected infinite graph and 
$G=\widetilde{G}/\Gamma$ a finite graph, where $\Gamma$ is a group of automorphisms of $\widetilde{G}$
such that the quotient map $\pi:\widetilde{G}\to G$
is a covering map, that is, 
it sends $N(x)$ bijectively onto $N(\pi(x))$ for every vertex $x$.
We prove that forcing commutes with taking preimages, so every closure in $G$
is the image of a closure in $\widetilde{G}$.
In the theory of covering spaces, an open set is evenly covered when 
its preimage is a disjoint union of copies of it. 
In the same spirit, we call a set $X\subseteq V(\widetilde{G})$ 
a \emph{sheet} if it is disjoint from, and has no edge to,
$\gamma X$ for every $\gamma\in \Gamma\setminus \{1\}$.
Then $\pi$ maps the subgraph induced by $X$ isomorphically onto 
the subgraph induced by $\pi(X)$, 
and $\pi^{-1}(\pi(X))=\bigsqcup_{\gamma\in\Gamma}\gamma X$.
The main tool of the paper (Theorem~\ref{thm:sheets}) says that
if the closure of every set of at most $s$ vertices of $\widetilde G$ 
is a sheet whenever it is connected, then no set of at most $s$ vertices 
is a zero forcing set of $G$, that is, $\Z(G)>s$.
Suppose $\widetilde{G}$ is periodic, that is, its vertices can be arranged
in columns $V_i$, where $V(\widetilde G)=\bigsqcup_{i\in\ZZ} V_i$,
and there is an automorphism $\sigma$ with $\sigma(V_i)=V_{i+1}$.
If $\Z(\widetilde G)$ is finite and $\widetilde G$ satisfies condition (W) below, then
$\Z(\widetilde G/\langle\tau\rangle)=\Z(\widetilde G)$
for every automorphism $\tau$ with $\tau(V_i)=V_{i+n}$ and $n$ large enough,
with an explicit bound on $n$ in terms of $h$, $g$ and $\Z(\widetilde G)$,
where the \emph{span} $h$ is the largest difference between the columns of the two ends of an edge,
and $g\ge h$ is the number of columns of a \emph{window}, a union of $g$ consecutive columns.

\begin{mainthm}[Theorem~\ref{thm:periodic_Z}]\label{thm:A}
  Let $\widetilde G$ be connected and periodic, with finite columns, and let $g\ge h$.
  Suppose that $\Z(\widetilde G)$ is finite and
  \begin{enumerate}
    \item[(W)] for every window $W$ of $g$ consecutive columns
      and every finite $X\subseteq V(\widetilde G)$ with $\cl(X)=X$,
      if $X\cap W\ne\emptyset$, then $\partial X\cap W\ne\emptyset$.
  \end{enumerate}
  If $\tau$ is an automorphism of $\widetilde G$ that moves every vertex $n$ columns to the right,
  where $n>2h$ and $n\ge g(\Z(\widetilde G)-2)+h+1$, then $\Z(\widetilde G/\langle\tau\rangle)=\Z(\widetilde G)$.
\end{mainthm}

The theorem says that beyond a certain threshold, every quotient $\widetilde G/\langle\tau\rangle$ has the same zero forcing number as $\widetilde G$.
The threshold on $n$ is explicit, so only finitely many graphs are left to check.
The proof needs no matrix construction.

We apply the theorem to three infinite families of graphs:
generalized Petersen graphs $\GP{n}{k}$, twisted cylinders $C_n\,\square_\varphi H$, and circulants $C_n(1,a)$.
We derive the zero forcing numbers of these graphs, some of which are known and some of which are new.

For the generalized Petersen graphs, we prove $\Z(\GP{n}{k})=2k+2$
for every $k\ge 1$ and $n\ge 2k^2+k+1$ (Theorem~\ref{thm:Z_finite}).
The upper bound was known~\cite{Alameda2018,RPT2020,Krishnan2026}.
The lower bound was known for $k=1$, where $\GP{n}{1}$ is a prism~\cite[Theorem~3.8]{AIM2008},
for $k=2$ and $n\ge 10$~\cite{RPT2020},
for $k=3$ and $n\ge 13$ and $k=4$ and $n\ge 18$~\cite{KMZ2026},
and for $k=5$ and $n\ge 23$~\cite{Boucher2026}.
It is new for all $k\ge 6$.
We also show that the least threshold $N_k$ satisfies $k^2-k+2\le N_k\le 2k^2+k+1$
(Proposition~\ref{prop:N_k}).
A computer search below the threshold gives $\Z(\GP{n}{k})$ for all $n$ when $k\le 5$
(Table~\ref{tab:values}).
This independently confirms the values of~\cite{KMZ2026} for $k=3,4$
and of~\cite{Boucher2026} for $k=5$.

For the twisted cylinders, we prove $\Z(C_n\,\square_\varphi H)=2m$
for every connected graph $H$ with $m$ vertices, every automorphism $\varphi$ of $H$,
and every $n\ge\max\{2m,3\}$ (Theorem~\ref{thm:cylinder}).
The threshold is best possible, since $\Z(C_n\,\square\,P_m)=\min\{n,2m\}$~\cite[Theorem~3.8]{AIM2008}.
Special cases were known, also for smaller $n$:
prisms and cylinders~\cite[Theorem~3.8]{AIM2008}, Möbius ladders~\cite[Proposition~3.9]{AIM2008},
tori~\cite[Theorem~2.6]{Benson2018} and $C_n\,\square\,K_m$~\cite[Theorem~4.1]{Duong2020}.
The result for general $H$ and $\varphi$ is new.

For the circulants, we prove $\Z(C_n(1,a))=2a$ for $a\ge 2$ and $n\ge 2a^2-a+1$
(Theorem~\ref{thm:circulant}),
and that the least such threshold $N'_a$ satisfies $a^2+1\le N'_a\le 2a^2-a+1$ for $a\ge 3$
(Proposition~\ref{prop:N_a}).
A computer search gives $N'_a=a^2+1$ for $a\le 6$ (Tables~\ref{tab:circ} and~\ref{tab:circ6}).
The case $a=2$ was known, since $C_n(1,2)$ is a consecutive circulant~\cite[Theorem~5.4]{DeaettMeyer2016},
and so was $\Z(C_{3m}(1,3))=6$ for $m>3$~\cite[Theorem~4.15]{Duong2020}.
The result is new for $a\ge 4$, and for $a=3$ when $n$ is not divisible by $3$.

\section{Closures and Boundaries}\label{sec:closures}
In this section we list some properties of closures for any graph $G$.
\begin{definition}
We say a set of vertices $X$ is \emph{stalled} if $\cl(X) = X$,
that is, if no vertex of $X$ has exactly one neighbor outside $X$.
\end{definition}
\begin{proposition}\label{prop:monotone}
  Let $S,T\subseteq V(G)$.
  If $S\subseteq T$, then $\cl(S)\subseteq\cl(T)$.
  Moreover, assuming every vertex of $G$ has finitely many neighbors, 
  $\cl(\cl(S))=\cl(S)$, that is, $\cl(S)$ is stalled.
\end{proposition}
\begin{proof}
  Let $B_0=S\subseteq B_1\subseteq\cdots$ and $B'_0=T\subseteq B'_1\subseteq\cdots$
  be the sequences defining $\cl(S)$ and $\cl(T)$.
  For the first statement we show $B_i\subseteq B'_i$ by induction on $i$.
  First, $B_0=S\subseteq T=B'_0$.
  Suppose $B_i\subseteq B'_i$, and let $w\in B_{i+1}\setminus B_i$,
  say $w$ is forced by $x\in B_i\subseteq B'_i$.
  Then every neighbor of $x$ other than $w$ lies in $B_i\subseteq B'_i$.
  So either $w\in B'_i$, or $w$ is the unique neighbor of $x$ outside $B'_i$, 
  in which case $w\in B'_{i+1}$.
  Hence $B_{i+1}\subseteq B'_{i+1}$.

  We prove the second statement by contradiction.
  Suppose $x\in\cl(S)$ has exactly one neighbor $w$ outside $\cl(S)$.
  The vertex $x$ and its finitely many neighbors other than $w$ lie in $\cl(S)$,
  so they all lie in some $B_j$.
  Then $w$ is the unique neighbor of $x$ outside $B_j$, 
  so $w\in B_{j+1}\subseteq\cl(S)$, a contradiction.
\end{proof}

\begin{remark}
  The second statement can fail in a graph with a vertex $x\in S$ of infinite degree:
  if all neighbors of $x$ except $w$ are forced in later and later rounds,
  then $x$ forces $w$ in no round, although $w$ is its only neighbor outside $\cl(S)$.
\end{remark}

\begin{definition}
  We say a set of vertices $X\subseteq V(G)$ is \emph{connected} 
  if any two vertices of $X$ are joined by a finite path in
  the subgraph induced by $X$.
  The \emph{connected components} of $X$ are its maximal connected subsets.
  We say two sets $X_1,X_2\subseteq V(G)$ \emph{touch} if some vertex of $X_1$
  is equal or adjacent to some vertex of $X_2$.
\end{definition}
Two distinct connected components of a set do not touch.
Conversely, if $X$ is the union of nonempty connected sets
and they do not touch pairwise,
then these sets are the connected components of $X$.

\begin{lemma}\label{lem:union_stalled}
  Let $\{X_\alpha\}$ be a family of stalled subsets of $V(G)$.
  If $X_\alpha$ and $X_\beta$ do not touch for all $\alpha \ne \beta$,
  then the union $U=\bigcup_\alpha X_\alpha$ is stalled. 
  (This is true for both finite and infinite families.)
\end{lemma}
\begin{proof}
  Suppose $U$ is not stalled.
  Then some $x\in U$ has exactly one neighbor $w$ outside $U$.
  Say $x\in X_\alpha$.
  Then $w$ is also outside $X_\alpha$, and all other neighbors of $x$ are in $U$.
  Since $X_\alpha$ touches no other $X_\beta$, 
  the vertex $x$ has no neighbor in such an $X_\beta$,
  so those other neighbors are in $X_\alpha$.
  Thus $w$ is the only neighbor of $x$ outside $X_\alpha$,
  and $X_\alpha$ is not stalled, a contradiction.
\end{proof}

\begin{lemma}\label{lem:closure_connected}
  If $S\subseteq V(G)$ is connected, then $\cl(S)$ is connected.
  Suppose $T, T'\subseteq V(G)$.
  If $\cl(T)$ and $\cl(T')$ are connected and touch, then $\cl(T\cup T')$ is also connected.
\end{lemma}
\begin{proof}
  Let $B_0=S\subseteq B_1\subseteq\cdots$ be as in the definition of $\cl(S)$.
  By induction, every $B_i$ is connected,
  since each vertex of $B_{i+1}\setminus B_i$ is adjacent to its forcer in $B_i$.
  Any two vertices of $\cl(S)$ lie in some $B_i$, so they are joined by a path in $B_i\subseteq\cl(S)$.
  Hence $\cl(S)$ is connected.

  For the second statement, let $S=\cl(T)\cup\cl(T')$,
  which is connected since $\cl(T)$ and $\cl(T')$ are connected and touch.
  By Proposition~\ref{prop:monotone}, $T\cup T'\subseteq S\subseteq\cl(T\cup T')$, thus
  \[
    \cl(T\cup T')\subseteq\cl(S)\subseteq\cl(\cl(T\cup T'))=\cl(T\cup T').
  \]
  Therefore, $\cl(T\cup T')=\cl(S)$, which is connected by the first statement.
\end{proof}

The next proposition says that the boundary never grows during the forcing process,
so it never has more than $|S|$ vertices.
For finite graphs, this bound is implicit in the count of
\emph{active} vertices in~\cite{Barioli2013}.
In periodic graphs, this bound limits how wide a connected closure can be
(Lemma~\ref{lem:width}).
\begin{proposition}\label{prop:boundary}
  Let $S\subseteq V(G)$
  and $S=B_0\subseteq B_1\subseteq B_2\subseteq\cdots$ be as in the definition of $\cl(S)$.
  Then $|\partial B_{i+1}|\le|\partial B_i|$ for every $i\ge 0$.
  In particular, $|\partial B_i|\le|S|$ for every $i\ge 0$, and $|\partial\cl(S)|\le|S|$.
\end{proposition}
\begin{proof}
  For each $w\in W_i=B_{i+1}\setminus B_i$, choose a forcer $f(w)\in B_i$.
  Then $f\colon W_i\to B_i$ is one-to-one, and $f(W_i)\subseteq\partial B_i$
  since $f(w)$ has the neighbor $w\notin B_i$.
  A vertex $v\in B_i\setminus\partial B_i$ is not in $\partial B_{i+1}$,
  since all neighbors of $v$ are in $B_i\subseteq B_{i+1}$.
  A vertex $f(w)$ is not in $\partial B_{i+1}$ either,
  since its only neighbor outside $B_i$ is $w\in B_{i+1}$.
  Therefore $\partial B_{i+1}\setminus W_i\subseteq\partial B_i\setminus f(W_i)$, and
  \[
    |\partial B_{i+1}|\le |W_i|+|\partial B_{i+1}\setminus W_i|
    \le |f(W_i)|+|\partial B_i\setminus f(W_i)|=|\partial B_i|.
  \]
  Since $|\partial B_0|=|\partial S|\le|S|$, by induction $|\partial B_i|\le|S|$ for every $i$.

  For the last inequality, let $D_i = \partial\cl(S)\cap B_i$. 
  If $v\in D_i$, then $v$ has a neighbor outside $\cl(S)$, 
  thus outside $B_{i}$, so $D_i\subseteq \partial B_i$ and $|D_i| \le |\partial B_i| \le |S|$.
  Then we have $\partial\cl(S) = \bigcup_{i} D_i$ and 
  $D_0 \subseteq D_1 \subseteq D_2 \subseteq \cdots$ with each $|D_i| \le |S|$.
  If $|S|$ is finite, then $D_i$ stabilize at some point, say for all $i\ge N$, $D_i=D_N$.
  In this case, $\partial \cl(S) = D_N$ and $|\partial \cl(S)| = |D_N| \le |S|$.
  If $|S|$ is infinite, then $\partial\cl(S)$ is a countable union of sets
  with at most $|S|$ elements, so $|\partial\cl(S)|\le\aleph_0\cdot|S|=|S|$.
\end{proof}

A \emph{two-way infinite path} in $G$ is a sequence $(p_j)_{j\in\ZZ}$ of distinct vertices
with $p_j\sim p_{j+1}$ for every $j$.

\begin{corollary}\label{cor:tracks}
  If $G$ contains $t$ pairwise disjoint two-way infinite paths, then $\Z(G)\ge 2t$.
\end{corollary}
\begin{proof}
  Let $S$ be a zero forcing set. 
  We show $|S|\ge 2t$, which is clear if $S$ is infinite,
  so we assume $S$ is finite.
  In this case each $B_i$ is finite by induction, since each vertex forces at most one new vertex.
  Choose two vertices on each of the $t$ paths.
  As $\cl(S)=V(G)$, those $2t$ vertices lie in $B_i$ for some $i$.
  Each path intersects $B_i$ at only finitely many vertices.
  Hence two vertices on each path at extreme ends are in $\partial B_i$.
  Since the paths are pairwise disjoint, this gives $2t$ distinct vertices of $\partial B_i$.
  Therefore, $2t\le|\partial B_i|\le|S|$ by Proposition~\ref{prop:boundary}.
\end{proof}

\section{Covers}\label{sec:covers}
In this section $\widetilde G$ is a graph in which every vertex has finitely many neighbors,
and $\Gamma$ is a group of automorphisms of $\widetilde G$
that moves every vertex to distance at least $3$.
In other words, for every vertex $x$ and every $\gamma\ne 1$ in $\Gamma$,
\begin{enumerate}
  \item[(i)] $\gamma x\ne x$,
  \item[(ii)] $\gamma x$ is not adjacent to $x$,
  \item[(iii)] $x$ and $\gamma x$ have no common neighbor.
\end{enumerate}
The \emph{quotient} $G=\widetilde G/\Gamma$ is the graph 
whose vertices are the orbits $\Gamma x$,
in which $\Gamma x \sim \Gamma y$ 
if some vertex of $\Gamma x$ is adjacent to some vertex of $\Gamma y$.
By (ii), no vertex is adjacent to a vertex of its own orbit, so $G$ has no loops.
Let $\pi\colon V(\widetilde G)\to V(G)$ be the map $\pi(x)=\Gamma x$.
By Lemma~\ref{lem:cover} below, $\pi$ is a covering map in the sense of topological graph theory~\cite{GrossTucker,Sunada2013}.

Automorphisms preserve adjacency,
so for all $\gamma \in \Gamma$, $u$ forces $v$ if and only if $\gamma u$ forces $\gamma v$.
Hence $\cl(\gamma T)=\gamma\cl(T)$ for every $T\subseteq V(\widetilde G)$,
and $\gamma X$ is stalled if and only if $X$ is.

\begin{lemma}\label{lem:cover}
  For every $x\in V(\widetilde G)$, the map $\pi$ sends $N(x)$ bijectively onto $N(\pi(x))$.
\end{lemma}
\begin{proof}
  If $y\in N(x)$, then $\pi(y)\ne\pi(x)$ by (ii), so $\pi(y)\in N(\pi(x))$.
  To show that $\pi\colon N(x) \to N(\pi(x))$ is surjective,
  let $\Gamma y\in N(\Gamma x)$. Then $\gamma_1 x\sim \gamma_2 y$ for some $\gamma_1,\gamma_2\in\Gamma$.
  Then $\gamma_1^{-1}\gamma_2 y\in N(x)$, and $\pi(\gamma_1^{-1}\gamma_2 y)=\Gamma y$.
  Next we prove injectivity.
  If $y_1,y_2\in N(x)$ and $\pi(y_1)=\pi(y_2)$, then $y_2=\gamma y_1$ for some $\gamma\in\Gamma$,
  and $x$ is a common neighbor of $y_1$ and $\gamma y_1$, so $\gamma=1$ by (iii) and $y_1=y_2$.
  Hence $\pi$ is a bijection from $N(x)$ onto $N(\pi(x))$.
\end{proof}

\begin{lemma}\label{lem:lift}
  For every $S\subseteq V(G)$, $\cl(\pi^{-1}(S))=\pi^{-1}(\cl(S))$.
  Consequently, $\cl(S)=\pi(\cl(\pi^{-1}(S)))$.
\end{lemma}
\begin{proof}
  Let $B_0=S\subseteq B_1\subseteq\cdots$ and $B'_0=\pi^{-1}(S)\subseteq B'_1\subseteq\cdots$
  be the sequences defining $\cl(S)$ and $\cl(\pi^{-1}(S))$, respectively.
  Once we show $B'_i=\pi^{-1}(B_i)$ by induction, taking the unions gives the first statement.
  The base case is true by definition.
  Suppose $B'_i=\pi^{-1}(B_i)$.
  Then for any $y\in V(\widetilde G)$, the bijection $\pi\colon N(y)\to N(\pi(y))$ of Lemma~\ref{lem:cover} restricts to a bijection
  $\pi\colon N(y)\setminus B'_i \to N(\pi(y))\setminus B_i$, since $B'_i=\pi^{-1}(B_i)$. Hence for $x\in N(y)$,
  we have 
  \begin{equation}\label{eq:neighbors}
      N(y)\setminus B'_i=\{x\} \iff N(\pi(y))\setminus B_i=\{\pi(x)\}.
  \end{equation}
  We claim, for any $x\in V(\widetilde G)$, 
  \begin{equation}\label{eq:forcing}
      x\in B'_{i+1}\setminus B'_i \iff \pi(x)\in B_{i+1}\setminus B_i.
  \end{equation}
  The forward direction is a direct consequence of (\ref{eq:neighbors}).
  For the other direction, we lift a forcer of $\pi(x)$ to a neighbor of $x$.
  Assume $\pi(x)\in B_{i+1}\setminus B_i$.
  Then there exists $z\in B_i$ such that $\pi(x)$ is the unique neighbor of $z$ outside $B_i$.
  Since $z\in N(\pi(x))$ and $\pi$ maps $N(x)$ onto $N(\pi(x))$,
  there exists $y\in N(x)$ such that $\pi(y)=z$, and $y\in B'_i$ since $z\in B_i$.
  Then the equivalence (\ref{eq:neighbors}) gives $N(y)\setminus B'_i=\{x\}$, so $x\in B'_{i+1}\setminus B'_i$.
  Thus the claim is true.
  Finally, we see that (\ref{eq:forcing}) gives $\pi^{-1}(B_{i+1}\setminus B_i) = B'_{i+1}\setminus B'_i$ and
  \[
  B'_{i+1} = B'_i\cup (B'_{i+1}\setminus B'_i) =\pi^{-1}(B_i)\cup\pi^{-1}(B_{i+1}\setminus B_i) = \pi^{-1}(B_{i+1}).
  \]
  The second statement follows since $\pi$ is surjective.
\end{proof}

\begin{definition}
  A set $X\subseteq V(\widetilde G)$ is a \emph{sheet} if
  $X$ and $\gamma X$ do not touch for every $\gamma\in\Gamma$ with $\gamma\ne 1$.
\end{definition}

\begin{example}\label{ex:path}
  Let $\widetilde G=P_\infty$ be the two-way infinite path with vertex set $\ZZ$ and edges $i\sim i+1$,
  and let $\Gamma=\langle\tau\rangle$, where $\tau(i)=i+n$ is the translation by $n$.
  Then (i)--(iii) hold if and only if $n\ge 3$, and $\widetilde G/\Gamma$ is the cycle $C_n$.
  For $n=2$, condition (iii) fails and $\pi$ does not map $N(i)$ bijectively onto $N(\pi(i))$.

  Let $n\ge 3$. A set $X$ of consecutive vertices is a sheet if and only if it does not touch $\tau X$,
  that is, if and only if it has at most $n-1$ vertices.
  For instance, $\cl(\{0\})=\{0\}$ is a sheet, but $\cl(\{0,1\})=\ZZ$ is not.
  This matches $\Z(C_n)=2$.
  In general, if $\Gamma\ne\{1\}$, the closure of a zero forcing set of $\widetilde G$,
  such as $\{0,1\}$ here, is $V(\widetilde G)$, which is not a sheet
  since $\gamma V(\widetilde G)=V(\widetilde G)$ for every $\gamma\in\Gamma$.
\end{example}

By Lemma~\ref{lem:lift}, we understand $\cl(S)$ by studying the connected components of $\cl(\pi^{-1}(S))$.
We build them by merging.
We start with one vertex of $\pi^{-1}(s)$ for each $s\in S$,
and whenever the closures of two of them touch, 
we can merge them after applying an automorphism.
Each merge lowers the number of sets, thus the process stops if $S$ is finite.

\begin{proposition}\label{prop:components}
  Let $S\subseteq V(G)$ be finite and nonempty, and suppose that $\cl(T)$ is a sheet
  for every $T\subseteq V(\widetilde G)$ such that $|T|\le|S|$ and $\cl(T)$ is connected.
  Then there are nonempty finite sets $T_1,\ldots,T_m\subseteq V(\widetilde G)$ with
  \[
    \pi(T_1)\cup\cdots\cup\pi(T_m)=S\qquad\text{and}\qquad |T_1|+\cdots+|T_m|=|S|
  \]
  such that the connected components of $\cl(\pi^{-1}(S))$ are the sets $\gamma\cl(T_r)$
  for $1\le r\le m$ and $\gamma\in\Gamma$.
  In particular, $\cl(S)=\pi(\cl(T_1))\cup\cdots\cup\pi(\cl(T_m))$.
\end{proposition}
\begin{proof}
  We obtain the sets $T_r$ by merging, keeping three conditions:
  \begin{enumerate}
    \item[(M1)] $\pi(T_1)\cup\cdots\cup\pi(T_m)=S$,
    \item[(M2)] $|T_1|+\cdots+|T_m|=|S|$,
    \item[(M3)] each $X_r=\cl(T_r)$ is connected.
  \end{enumerate}
  We start with $m=|S|$ and singletons $T_r=\{x_r\}$ by choosing $x_r\in \pi^{-1}(v_r)$
  for each $v_r\in S$.
  (M1) and (M2) hold trivially, and (M3) follows from Lemma~\ref{lem:closure_connected}.

  While (M1)--(M3) hold, $|T_r|\le|S|$ by (M2), so $X_r$ is a sheet by hypothesis,
  and the sets $\gamma X_r$ with $\gamma\in\Gamma$ pairwise do not touch.
  Now for $r\ne r'$, if $\gamma_1 X_r$ touches $\gamma_2 X_{r'}$,
  then $X_r$ touches $\lambda X_{r'}=\cl(\lambda T_{r'})$ with $\lambda=\gamma_1^{-1}\gamma_2$.
  In this case, we replace $T_r$ and $T_{r'}$ by $T_r\cup\lambda T_{r'}$.
  Its closure is connected by Lemma~\ref{lem:closure_connected},
  and its image is $\pi(T_r)\cup\pi(T_{r'})$.
  By (M1) and (M2), 
  \[
  |S|\le|\pi(T_1)|+\cdots+|\pi(T_m)|\le|T_1|+\cdots+|T_m|=|S|,
  \]
  so the sets $\pi(T_r)$ are pairwise disjoint.
  Since $\pi(\lambda T_{r'})=\pi(T_{r'})$, the sets $T_r$ and $\lambda T_{r'}$ are disjoint,
  so their union has $|T_r|+|T_{r'}|$ elements, and (M1)--(M3) still hold, with one set fewer.
  Hence after fewer than $|S|$ merges, no two of the sets $\gamma X_r$ touch.

  The following set $U$ has connected components $\gamma X_r$ by construction, 
  and we prove that $\cl(\pi^{-1}(S))=U$.
  \[
    U=\bigcup_{r=1}^m\bigcup_{\gamma\in\Gamma}\gamma X_r.
  \]
  First, each $\gamma X_r=\cl(\gamma T_r)$ is stalled by Proposition~\ref{prop:monotone},
  so $U$ is stalled by Lemma~\ref{lem:union_stalled}.
  Second, two vertices with the same image lie in the same orbit,
  so by (M1) every vertex of $\pi^{-1}(S)$ has the form $\gamma x$ with $x\in T_r$ for some $r$.
  Hence $\pi^{-1}(S)\subseteq U$, and
  $\cl(\pi^{-1}(S))\subseteq\cl(U)=U$ by Proposition~\ref{prop:monotone}.
  Conversely, for all $r$ and $\gamma$, $\gamma T_r\subseteq\pi^{-1}(S)$,
  so $\gamma X_r=\cl(\gamma T_r)\subseteq\cl(\pi^{-1}(S))$.
  Hence $\cl(\pi^{-1}(S))=U$.
  Finally, $\cl(S)=\pi(\cl(\pi^{-1}(S)))=\pi(U)=\pi(X_1)\cup\cdots\cup\pi(X_m)$ by Lemma~\ref{lem:lift}.
\end{proof}

\begin{theorem}\label{thm:sheets}
  Suppose that $\widetilde G$ is connected and $\Gamma\ne\{1\}$.
  Let $s\ge 1$, and suppose that $\cl(T)$ is a sheet
  for every $T\subseteq V(\widetilde G)$ such that $|T|\le s$ and $\cl(T)$ is connected.
  Then $\Z(G)\ge s+1$.
\end{theorem}
\begin{proof}
  Suppose $S$ is a zero forcing set of $G$ with $|S|\le s$. Then $S\ne\emptyset$, and
  $\cl(\pi^{-1}(S))=\pi^{-1}(V(G))=V(\widetilde G)$ by Lemma~\ref{lem:lift}.
  Since $\widetilde G$ is connected and $|S|\le s$, Proposition~\ref{prop:components} gives $V(\widetilde G)=\cl(T_1)$ with $|T_1|\le s$.
  So $V(\widetilde G)$ is a sheet by hypothesis,
  which is impossible since $\gamma V(\widetilde G)=V(\widetilde G)$ for $\gamma\ne 1$ in $\Gamma$.
\end{proof}

\section{Periodic Covers}\label{sec:periodic}

In this section we keep the setting of Section~\ref{sec:covers}
and give conditions on $\widetilde G$ that make closures of small sets sheets,
when $\widetilde G$ is periodic.

\begin{definition}
  We say $\widetilde G$ is \emph{periodic} if there are
  a partition $V(\widetilde G)=\bigcup_{i\in\ZZ}V_i$
  and an automorphism $\sigma$ of $\widetilde G$ with $\sigma(V_i)=V_{i+1}$ for every $i$.
  We call $V_i$ the $i$th \emph{column},
  write $\col(x)=i$ if $x\in V_i$, and say that a set $X$ \emph{meets} column $i$
  if $X\cap V_i\ne\emptyset$.
  The \emph{span} of $\widetilde G$ is $\sup\{\lvert\col(x)-\col(y)\rvert : xy\in E(\widetilde G)\}$,
  the largest number of columns that a single edge can cross.
  For $c\le d$ in $\ZZ\cup\{\pm\infty\}$, let $V[c,d]=\{x\in V(\widetilde G):c\le\col(x)\le d\}$,
  written $V[c,\infty)$ if $d=\infty$ and $V(-\infty,d]$ if $c=-\infty$.
  A \emph{window} of $g$ columns is a set $V[i,i+g-1]$.
  For a finite nonempty $X\subseteq V(\widetilde G)$, the \emph{width} of $X$ is
  \[
    \wid(X)=\max_{x\in X}\col(x)-\min_{x\in X}\col(x)+1.
  \]
  An automorphism $\tau$ with $\tau(V_i)=V_{i+n}$ for every $i$ is called a \emph{shift by $n$}.
\end{definition}

\begin{figure}[ht]
  \centering
  \begin{tikzpicture}[x=9mm, y=9mm,
      vertex/.style = {circle, fill=white, draw=black, line width=.5pt,
                       inner sep=0pt, minimum size=1.8mm},
      edge/.style   = {line width=.5pt, black!60}]
    \fill[black!10, rounded corners=3pt] (-0.3,-0.3) rectangle (0.3,1.3);
    \foreach \y in {0,1} {
      \draw[edge] (-1.55,\y) -- (5.55,\y);
      \node at (-1.95,\y) {$\cdots$};
      \node at (5.95,\y) {$\cdots$};
    }
    \foreach \i in {-1,...,5} {
      \draw[edge] (\i,0) -- (\i,1);
      \node[vertex] at (\i,0) {};
      \node[vertex] at (\i,1) {};
      \node[font=\small, black!60] at (\i,-1.0) {$\i$};
    }
    \node[vertex, fill=black] at (0,0) {};
    \node[vertex, fill=black] at (1,0) {};
    \node[vertex, fill=black] at (4,1) {};
    \node[font=\small] at (0,-0.5) {$x$};
    \node[font=\small] at (1,-0.5) {$\sigma(x)$};
    \node[font=\small] at (4,1.5) {$\tau(x)$};
    \node[font=\small] at (0,1.5) {$V_0$};
  \end{tikzpicture}
  \caption{The infinite ladder $P_\infty\,\square\,K_2$ with columns indexed by $\ZZ$ and the column $V_0$ shaded.
    The automorphism $\sigma$ maps $x$ to $\sigma(x)$,
    and the shift $\tau$ by $n=4$ maps $x$ to $\tau(x)$ in the other row.}
  \label{fig:ladder}
\end{figure}
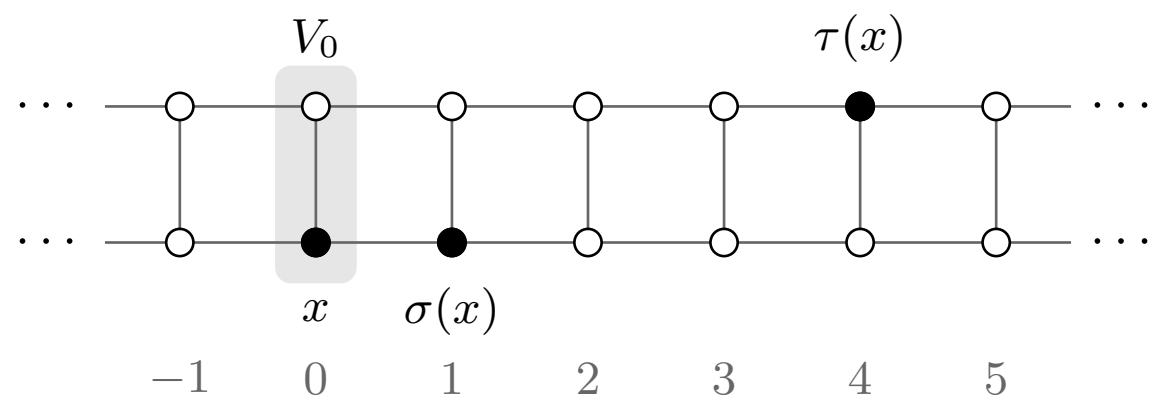

\begin{example}\label{ex:ladder}
  The \emph{infinite ladder} $P_\infty\,\square\,K_2$ has vertex set $\ZZ\times\{0,1\}$, with
  $(i,0)\sim(i,1)$ and $(i,a)\sim(i+1,a)$ for all $i\in\ZZ$ and $a\in\{0,1\}$ (Figure~\ref{fig:ladder}).
  It is periodic with columns $V_i=\{(i,0),(i,1)\}$, the automorphism $\sigma(i,a)=(i+1,a)$, and span $1$.
  Every $\sigma^n$ is a shift by $n$, but so is $\tau(i,a)=(i+n,1-a)$, which also switches the two rows.
  For $n\ge 3$, the quotient by $\langle\sigma^n\rangle$ is the prism $C_n\,\square\,K_2$,
  and the quotient by $\langle\tau\rangle$ is the Möbius ladder on $2n$ vertices.
  Section~\ref{sec:cylinders} treats $P_\infty\,\square\,H$ for every finite connected graph $H$ in the same way.
\end{example}

For the rest of the section, we assume $\widetilde G$ is periodic with finite span $h$.
(The span is finite if $|V_i|<\infty$.)
Since the columns $V_i$ are disjoint and nonempty, $\widetilde G$ is infinite.
Suppose $\tau$ is a shift by $n>2h$, and let $\Gamma=\langle\tau\rangle$.
For $j\ne 0$, the automorphism $\tau^j$ changes the column by $\lvert j\rvert n>2h$,
while a path with at most two edges changes it by at most $2h$.
Hence $\Gamma$ acts freely and moves every vertex to distance at least $3$,
so the results of Section~\ref{sec:covers} apply to $\tilde G/\Gamma$.

\begin{lemma}\label{lem:narrow}
  Let $\tau$ be a shift by $n$.
  Every finite nonempty $X\subseteq V(\widetilde G)$ with $\wid(X)\le n-h$
  is a sheet for $\langle\tau\rangle$.
\end{lemma}
\begin{proof}
  Let $a$ and $b$ be the leftmost and rightmost columns that $X$ meets, so $b-a+1=\wid(X)$.
  For $j\ge 1$, $\tau^jX$ meets only columns $a+jn,\ldots,b+jn$, and
  \[
    (a+jn)-b\ge (a+n)-b=n-\wid(X)+1\ge h+1.
  \]
  Since an edge crosses at most $h$ columns,
  no vertex of $\tau^jX$ is equal or adjacent to a vertex of $X$, so $X$ and $\tau^jX$ do not touch.
  The proof is similar for $j\le -1$.
\end{proof}

\begin{lemma}\label{lem:dichotomy}
  Suppose $\widetilde G$ is connected, every column of $\widetilde G$ is finite, and $\Z(\widetilde G)$ is finite.
  Then for every finite $T\subseteq V(\widetilde G)$, $\cl(T)$ is either finite or all of $V(\widetilde G)$.
\end{lemma}
\begin{proof}
  Since $V_0\cup V_1$ is finite and $\widetilde G$ is connected, there is an $r\ge 1$
  such that any two vertices of $V_0\cup V_1$ are joined by a path in $V[-r,r]$.
  Since $\widetilde G$ is periodic, any two vertices in $V_i\cup V_{i+1}$ are joined by a path in $V[i-r,i+r]$.
  Thus, for $a\le b$, any two vertices of $V[a,b]$ are joined by a path in $V[a-r,b+r]$.
  Such a path can be constructed by picking intermediate vertices in
  each column and concatenating paths between them.

  Suppose $X=\cl(T)$ is infinite. By Proposition~\ref{prop:boundary}, $\partial X$ is finite.
  Thus there exists $M>0$ such that $\partial X\subseteq V[-M+1,M-1]$.
  Then either $X\cap V[M,\infty)$ or $X\cap V(-\infty,-M]$ is infinite.
  By symmetry, we may assume the former and pick a vertex $x\in X\cap V[M+r,\infty)$.
  For any $y\in V[M+r,\infty)$, there is a path $P$ in $V[M,\infty)$ joining $x$ and $y$.
  Since every vertex in $X\cap V[M,\infty)$ is outside $\partial X$, it has all its neighbors in $X$.
  Consequently, every vertex in the path $P$ lies in $X$ by induction starting from $x$, 
  in particular, $y\in X$. 
  Therefore, $V[M+r,\infty)\subseteq X$.

  Since $\Z(\widetilde G)$ is finite, there is a finite zero forcing set $S$,
  and some power $\sigma^k$ maps $S$ into $V[M+r,\infty)$,
  where $\sigma^k(S)$ is again a zero forcing set.
  Hence $X$ contains a zero forcing set, and since $X$ is stalled, $X=V(\widetilde G)$.
\end{proof}

\begin{definition}
  Let $g\ge 1$. We say $\widetilde G$ has the \emph{$g$-window property}
  if for every finite stalled $X\subseteq V(\widetilde G)$
  and every window $W$ of $g$ columns,
  if $X\cap W\ne\emptyset$, then $\partial X\cap W\ne\emptyset$.
\end{definition}

The $g$-window property implies the $(g+1)$-window property.
We call the least such $g$ the \emph{window number} of $\widetilde G$.
The window number determines the threshold on $n$ in Theorem~\ref{thm:periodic_Z} below.
Since Lemma~\ref{lem:width} below needs $g\ge h$,
we may take $g$ to be the larger of the window number and $h$.
In each of our applications in Sections~\ref{sec:gp}, \ref{sec:cylinders} and~\ref{sec:circulants}
the window number is at most $h$, so $g=h$.
There this follows from the next lemma.
For $W\subseteq V(\widetilde G)$, let $N[W]$ be the set of vertices in $W$ or adjacent to a vertex of $W$.

\begin{lemma}\label{lem:window_criterion}
  Let $g\ge 1$, and let $W$ be a window of $g$ columns.
  Suppose that $W$ is connected and $N[W]$ is a zero forcing set of $\widetilde G$.
  Then $\widetilde G$ has the $g$-window property.
\end{lemma}
\begin{proof}
  Every window of $g$ columns is $\sigma^j(W)$ for some $j\in\ZZ$,
  and they all have the same property as $W$.
  So it suffices to prove the window property for $W$.
  Let $X$ be finite and stalled, and suppose $X\cap W\ne\emptyset$ but $\partial X\cap W=\emptyset$.
  If $x\in X\cap W$, then $x\notin\partial X$, so every neighbor of $x$ is in $X$.
  Since $W$ is connected, repeating this along paths in $W$ gives $W\subseteq X$.
  Then every vertex of $W$ is outside $\partial X$, so $N[W]\subseteq X$.
  Thus $X=\cl(S)\supseteq\cl(N[W])=V(\widetilde G)$, which contradicts the finiteness of $X$.
\end{proof}

\begin{lemma}\label{lem:width}
  Suppose $\widetilde G$ has the $g$-window property for some $g\ge h$.
  If $X\subseteq V(\widetilde G)$ is finite, nonempty, stalled and connected, then
  \[
    \wid(X)\le g(|\partial X|-1)+1.
  \]
\end{lemma}
\begin{proof}
  Let $L$ and $R$ be the least and greatest columns met by $X$, so $\wid(X)=R-L+1$.
  Let $W=V[L-g+1,L]$, a window of $g$ columns meeting $X$ only in column $L$. 
  By the $g$-window property, $\partial X$ meets column $L$.
  Similarly, $\partial X$ meets column $R$.
  Let $d_1<\cdots<d_m$ be the columns met by $\partial X$.
  Then $d_1=L$, $d_m=R$, and $m\le|\partial X|$.

  Suppose $d_{i+1}-d_i>g$ for some $i$.
  Then we have a window $W=V[d_i+1,d_i+g]$ of $g$ columns such that $\partial X\cap W=\emptyset$,
  which implies $X\cap W=\emptyset$ by the $g$-window property.
  Then $X$ is broken into two parts separated by the interval $W$. 
  The two parts, each of which is nonempty since $X$ meets columns $d_i$ and $d_{i+1}$, 
  cannot be joined by an edge since an edge can cross at most $h\le g$ columns.
  This is a contradiction to the connectedness of $X$.
  Hence $d_{i+1}-d_i\le g$ for each $i<m$, and
  \[
    \wid(X)-1=d_m-d_1=\sum_{i<m}(d_{i+1}-d_i)\le g(m-1)\le g(|\partial X|-1). \qedhere
  \]
\end{proof}

\begin{theorem}\label{thm:periodic_Z}
  Let $\widetilde G$ be a periodic graph with span $h$ such that
  \begin{enumerate}
    \item[(a)] $\widetilde G$ is connected and its columns are finite,
    \item[(b)] $\widetilde G$ has the $g$-window property for some $g\ge h$, and
    \item[(c)] $\Z(\widetilde G)$ is finite.
  \end{enumerate}
  If $\tau$ is a shift by $n$ with $n>2h$ and $n\ge g(\Z(\widetilde G)-2)+h+1$, then
  \[
    \Z(G)=\Z(\widetilde G).
  \]
  where $G=\widetilde G/\langle\tau\rangle$.
\end{theorem}
\begin{proof}
  Since $n>2h$, the results of Section~\ref{sec:covers} apply to $\Gamma=\langle\tau\rangle$,
  and $\Gamma\ne \{1\}$ since $\tau$ moves columns nontrivially.
  First we prove $\Z(G)\ge \Z(\widetilde G)$.
  For this lower bound, we may assume $\Z(\widetilde G)\ge 2$, otherwise the lower bound is trivial.
  We let $s=\Z(\widetilde G)-1\ge 1$ and verify the conditions of Theorem~\ref{thm:sheets}.
  So we assume $T\subseteq V(\widetilde G)$, $|T|\le s$ and $X=\cl(T)$ is connected, 
  and show that $X$ is a sheet.
  If $T=\emptyset$, then $X=\emptyset$ is a sheet, so assume $T\ne\emptyset$.
  Since $|T|<\Z(\widetilde G)$, we have $X\ne V(\widetilde G)$, thus $X$ is finite by Lemma~\ref{lem:dichotomy}.
  Moreover, $X$ is nonempty, stalled, and $|\partial X|\le|T|\le s$ by Proposition~\ref{prop:boundary}.
  Also by Lemma~\ref{lem:width}, 
  \[
    \wid(X)\le g(|\partial X|-1)+1\le g(s-1)+1= g(\Z(\widetilde G)-2)+1\le n-h.
  \]
  Hence $X$ is a sheet by Lemma~\ref{lem:narrow}.
  Since $\widetilde G$ is connected, Theorem~\ref{thm:sheets} gives $\Z(G)\ge s+1=\Z(\widetilde G)$.

  For the upper bound, let $S$ be a zero forcing set of $\widetilde G$ with $|S|=\Z(\widetilde G)$.
  Since $S\subseteq\pi^{-1}(\pi(S))$, Lemma~\ref{lem:lift} and Proposition~\ref{prop:monotone} give
  \[
    \cl(\pi(S))=\pi(\cl(\pi^{-1}(\pi(S))))\supseteq\pi(\cl(S))=\pi(V(\widetilde G))=V(G).
  \]
  Thus $\pi(S)$ is a zero forcing set of $G$, and $\Z(G)\le|\pi(S)|\le|S|=\Z(\widetilde G)$.
\end{proof}

In the remaining sections we determine $\Z(\widetilde G)$ and $g$ for three families of covers.

\section{Generalized Petersen Graphs}\label{sec:gp}

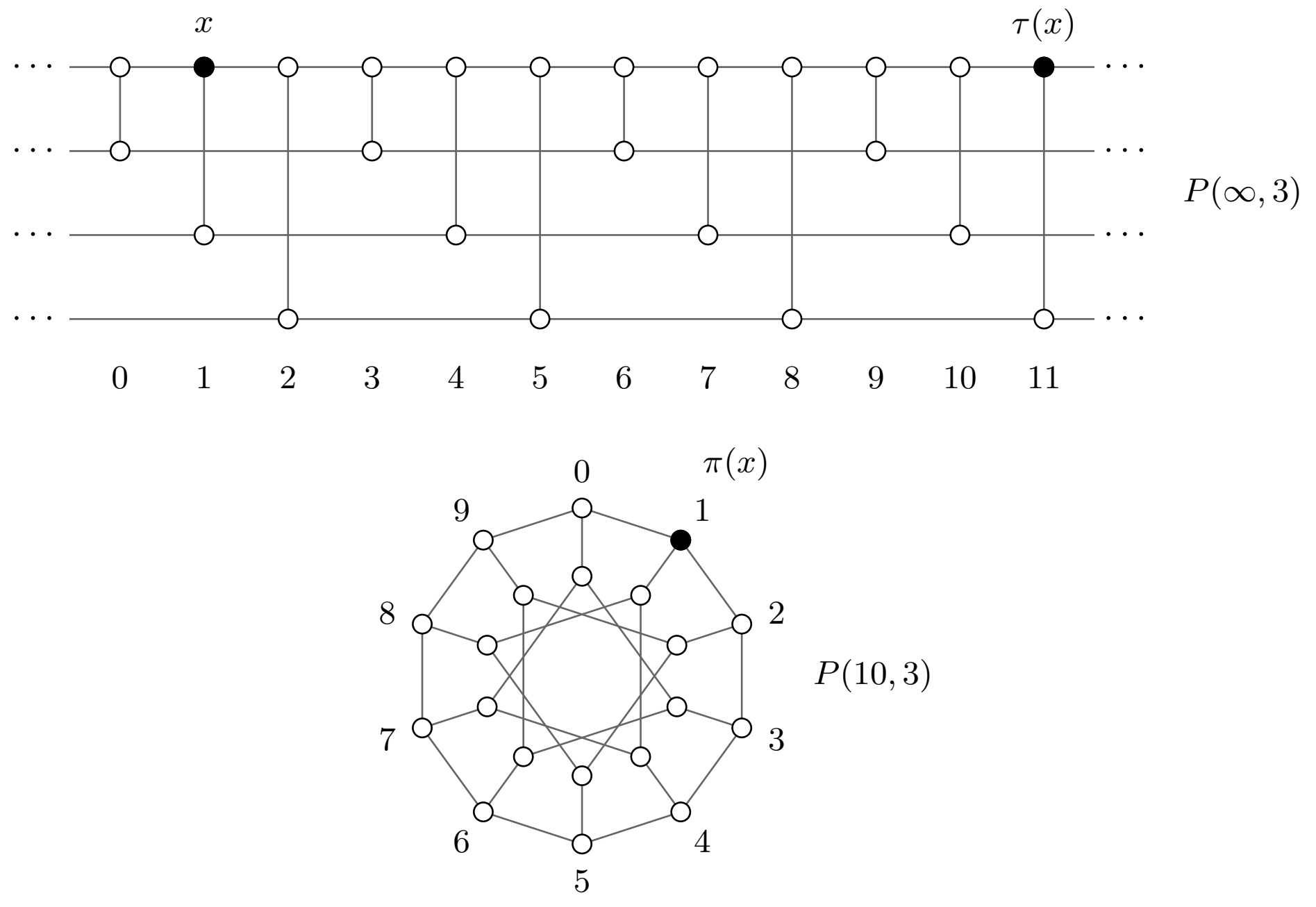
\begin{figure}[ht]
  \centering
  \begin{tikzpicture}[x=8mm, y=8mm,
      vertex/.style = {circle, fill=white, draw=black, line width=.5pt,
                       inner sep=0pt, minimum size=1.8mm},
      edge/.style   = {line width=.5pt, black!60}]
    % top: the cover P(infty,3)
    % u_i is drawn in column i and row 3, and v_i in column i and row 2 - (i mod 3)
    \begin{scope}
      \begin{scope}
        \clip (-0.6,-0.5) rectangle (11.6,3.5);
        \foreach \y in {0,1,2,3} \draw[edge] (-1,\y) -- (12,\y);
      \end{scope}
      \foreach \i in {0,...,11} {
        \pgfmathtruncatemacro{\r}{2-mod(\i,3)}
        \draw[edge] (\i,3) -- (\i,\r);
      }
      \foreach \y in {0,1,2,3} {
        \node at (-1,\y) {$\cdots$};
        \node at (12,\y) {$\cdots$};
      }
      \foreach \i in {0,...,11} {
        \pgfmathtruncatemacro{\r}{2-mod(\i,3)}
        \node[vertex] at (\i,3) {};
        \node[vertex] at (\i,\r) {};
        \node[font=\small] at (\i,-0.7) {$\i$};
      }
      \node[vertex, fill=black] at (1,3) {};
      \node[vertex, fill=black] at (11,3) {};
      \node[font=\small] at (1,3.5) {$x$};
      \node[font=\small] at (11,3.5) {$\tau(x)$};
      \node[font=\small, anchor=west] at (12.5,1.5) {$\Pinf{3}$};
    \end{scope}
    % bottom: the generalized Petersen graph P(10,3)
    \begin{scope}[yshift=-3.4cm, xshift=4.4cm]
      \foreach \i in {0,...,9} {
        \coordinate (u\i) at ({90-36*\i}:1.6cm);
        \coordinate (v\i) at ({90-36*\i}:0.95cm);
      }
      \foreach \i in {0,...,9} {
        \pgfmathtruncatemacro{\j}{mod(\i+1,10)}
        \pgfmathtruncatemacro{\l}{mod(\i+3,10)}
        \draw[edge] (u\i) -- (u\j);
        \draw[edge] (u\i) -- (v\i);
        \draw[edge] (v\i) -- (v\l);
      }
      \foreach \i in {0,...,9} {
        \node[vertex] at (u\i) {};
        \node[vertex] at (v\i) {};
        \node[font=\small] at ({90-36*\i}:1.95cm) {$\i$};
      }
      \node[vertex, fill=black] at (u1) {};
      \node[font=\small] at (54:2.5cm) {$\pi(x)$};
      \node[font=\small, anchor=west] at (2.6,0) {$\GP{10}{3}$};
    \end{scope}
  \end{tikzpicture}
  \caption{The cover $\Pinf{3}$ (top) and the generalized Petersen graph $\GP{10}{3}$ (bottom).}
  \label{fig:gp}
\end{figure}

Fix $k\ge 1$.
For $n>2k$, the \emph{generalized Petersen graph} $\GP{n}{k}$
has $n$ outer vertices $u_0,\ldots, u_{n-1}$, $n$ inner vertices $v_0,\ldots, v_{n-1}$, and edges
\[
  u_i\sim u_{i+1}\ \text{(outer edges)},\quad
  u_i\sim v_i\ \text{(spokes)},\quad
  v_i\sim v_{i+k}\ \text{(inner edges)},
\]
indices modulo~$n$.
Its cover, the \emph{infinite generalized Petersen graph} $\Pinf{k}$,
has vertex set $\{u_i: i\in\ZZ\}\cup\{v_i: i\in\ZZ\}$ and the same adjacencies.
The graph $\Pinf{k}$ is $3$-regular, with
\[
  N(u_i)=\{u_{i-1},u_{i+1},v_i\}, \qquad N(v_i)=\{v_{i-k},v_{i+k},u_i\}.
\]
The outer vertices induce a two-way infinite path,
and the inner vertices induce $k$ of them, one for each residue of the index modulo $k$.
With the columns $V_i=\{u_i,v_i\}$ and the translation $\sigma(u_i)=u_{i+1}$, $\sigma(v_i)=v_{i+1}$,
the graph $\Pinf{k}$ is periodic with span $h=k$.
The map $\tau=\sigma^n$ is an automorphism of $\Pinf{k}$ and a shift by $n$,
and $\pi(u_i)=u_{i\bmod n}$, $\pi(v_i)=v_{i\bmod n}$ identifies $\Pinf{k}/\langle\tau\rangle$ with $\GP{n}{k}$.
See Figure~\ref{fig:gp}.

The next lemma is the analogue for $\Pinf{k}$ of the bound $\Z(\GP{n}{k})\le 2k+2$
proved by Alameda et al.~\cite[Theorem~2.2]{Alameda2018} and by Rashidi et al.~\cite[Theorem~2.2]{RPT2020},
and the proof below is essentially the same.

\begin{lemma}\label{lem:outer_forced}
  Any $2k+2$ consecutive outer vertices form a zero forcing set of $\Pinf{k}$.
\end{lemma}

\begin{proof}
  By translation, we may assume $u_1,\ldots,u_{2k+2}$ are black.
  For $2\le i\le 2k+1$, $u_i$ forces $v_i$, since $u_{i-1}$ and $u_{i+1}$ are black.
  Then $v_{k+2}$ forces $v_{2k+2}$, since $u_{k+2}$ and $v_2$ are black,
  and $u_{2k+2}$ forces $u_{2k+3}$, since $u_{2k+1}$ and $v_{2k+2}$ are black.
  Symmetrically, $v_{k+1}$ forces $v_1$, and $u_1$ forces $u_0$.
  Repeating this on both sides, every vertex becomes black.
\end{proof}

\begin{lemma}\label{lem:window}
  The graph $\Pinf{k}$ has the $k$-window property.
\end{lemma}

\begin{proof}
  We apply Lemma~\ref{lem:window_criterion} to $W=V[1,k]$,
  which is connected by the outer path $u_1,\ldots,u_k$ and the spokes.
  Its closed neighborhood is
  \[
    N[W]=\{u_i: 0\le i\le k+1\}\cup\{v_i: 1-k\le i\le 2k\}.
  \]
  For $j=k+1,\ldots,2k$ in turn, $u_j$ forces $u_{j+1}$,
  since its other neighbors $u_{j-1}$ and $v_j\in N[W]$ are black.
  Hence $\cl(N[W])$ contains the set $\{u_0,\ldots,u_{2k+1}\}$,
  which is a zero forcing set by Lemma~\ref{lem:outer_forced},
  so $\cl(N[W])=\cl(\cl(N[W]))=V(\Pinf{k})$ by Proposition~\ref{prop:monotone}.
\end{proof}

\begin{proposition}\label{prop:Z_inf}
  $\Z(\Pinf{k})=2k+2$.
\end{proposition}

\begin{proof}
  By Lemma~\ref{lem:outer_forced}, $\Z(\Pinf{k})\le 2k+2$.
  The outer path and the $k$ inner paths are $k+1$ disjoint two-way infinite paths,
  hence $\Z(\Pinf{k})\ge 2k+2$ by Corollary~\ref{cor:tracks}.
\end{proof}

\begin{theorem}\label{thm:Z_finite}
  Let $n\ge 2k^2+k+1$. Then
  \[
    \Z(\GP{n}{k})=2k+2.
  \]
\end{theorem}

\begin{proof}
  We apply Theorem~\ref{thm:periodic_Z} to $\widetilde G=\Pinf{k}$ with $h=g=k$.
  Condition~(a) holds since $\Pinf{k}$ is connected and each column has two vertices,
  (b) holds by Lemma~\ref{lem:window},
  and (c) holds since $\Z(\Pinf{k})=2k+2$ by Proposition~\ref{prop:Z_inf}.
  The threshold in Theorem~\ref{thm:periodic_Z} is
  $g(\Z(\Pinf{k})-2)+h+1=k(2k+2-2)+k+1=2k^2+k+1$.
  Since $n>2k$ and $n\ge 2k^2+k+1$, Theorem~\ref{thm:periodic_Z} gives $\Z(\GP{n}{k})=2k+2$.
  The upper bound $\Z(\GP{n}{k})\le 2k+2$ is also \cite[Theorem~2.2]{Alameda2018}.
\end{proof}

As a consequence of the theorem, for each $k\ge 1$ there is a least integer $N_k$
such that $\Z(\GP{n}{k})=2k+2$ for all $n\ge N_k$.

\begin{proposition}\label{prop:N_k}
  For every $k\ge 1$,
  \[
    k^2-k+2\le N_k\le 2k^2+k+1.
  \]
\end{proposition}

\begin{proof}
  The upper bound is Theorem~\ref{thm:Z_finite}.
  For $k\le 2$ the lower bound follows from Table~\ref{tab:values}, since $N_1=4$ and $N_2=9$.
  For $k\ge 3$, let $n=k^2-k+1$.
  Then $n>2k$, so $\GP{n}{k}$ and $\GP{n}{k-1}$ are defined, and $k(k-1)=n-1\equiv -1\pmod n$.
  By Watkins' theorem~\cite{Watkins1969}, $\GP{n}{k}\cong\GP{n}{l}$ whenever $kl\equiv\pm 1\pmod n$,
  so $\GP{n}{k}\cong\GP{n}{k-1}$.
  Therefore, by \cite[Theorem~2.2]{Alameda2018},
  \[
    \Z(\GP{n}{k})=\Z(\GP{n}{k-1})\le 2k,
  \]
  proving $N_k>n$.
\end{proof}

For small $k$ we determine $\Z(\GP{n}{k})$ for the remaining $n$, that is, for $2k<n<2k^2+k+1$,
by an exhaustive computer search.
For $k=1$ this is only the prism $\GP{3}{1}=C_3\square K_2$, and $\Z(\GP{3}{1})=3$ by \cite[Theorem~3.8]{AIM2008}.
For $2\le k\le 5$ the results are in Table~\ref{tab:values}.
For $k=3$ and $k=4$ these values agree with those of
Krishnan, Ma and Zhang~\cite[Theorem~1.1]{KMZ2026},
and for $k=5$ with those of Boucher~\cite{Boucher2026}.
Rashidi et al.~\cite{RPT2020} stated $\Z(\GP{n}{3})=8$ for $n\ge 12$ with an incomplete proof,
and Krishnan~\cite{Krishnan2026} showed $\Z(\GP{12}{3})=7$, as in Table~\ref{tab:values}.
Bjorkman et al.~\cite[Table~1, $\ell=0$]{Bjorkman2025} list $\Z(\GP{n}{2})$ for $5\le n\le 9$
and $\Z(\GP{n}{3})$ for $7\le n\le 11$.
These agree with Table~\ref{tab:values} except $\Z(\GP{7}{2})$, which they give as $5$.
In fact $\Z(\GP{7}{2})=6$, since $\GP{7}{2}\cong\GP{7}{3}$ by Watkins' theorem
and their table gives $\Z(\GP{7}{3})=6$.
For $k=6$ we have $32\le N_6\le 79$ by Proposition~\ref{prop:N_k}.

\begin{table}[ht]
  \caption{The zero forcing number $\Z(\GP{n}{k})$ for $k\le 5$.}
  \label{tab:values}
  \setlength{\tabcolsep}{2.7pt}
  \begin{tabular}{c|*{20}{c}|c|c}
    $k\backslash n$ & 3 & 4 & 5 & 6 & 7 & 8 & 9 & 10 & 11 & 12 & 13 & 14 & 15 & 16 & 17 & 18 & 19 & 20 & 21 & 22 & $\ge 23$ & $N_k$ \\
    \hline
    1 & 3 & 4 & 4 & 4 & 4 & 4 & 4 & 4 & 4 & 4 & 4 & 4 & 4 & 4 & 4 & 4 & 4 & 4 & 4 & 4 & 4 & 4 \\
    2 &   &   & 5 & 4 & 6 & 5 & 6 & 6 & 6 & 6 & 6 & 6 & 6 & 6 & 6 & 6 & 6 & 6 & 6 & 6 & 6 & 9 \\
    3 &   &   &   &   & 6 & 6 & 6 & 8 & 7 & 7 & 8 & 8 & 8 & 8 & 8 & 8 & 8 & 8 & 8 & 8 & 8 & 13 \\
    4 &   &   &   &   &   &   & 6 & 6 & 7 & 6 & 8 & 8 & 9 & 8 & 9 & 10 & 10 & 10 & 10 & 10 & 10 & 18 \\
    5 &   &   &   &   &   &   &   &   & 6 & 8 & 9 & 8 & 8 & 8 & 10 & 10 & 10 & 10 & 10 & 11 & 12 & 23 \\
  \end{tabular}
\end{table}

\section{Cylinders and Twisted Cylinders}\label{sec:cylinders}

Let $H$ be a finite connected graph with $m$ vertices.
The graph $\widetilde G=P_\infty\,\square\, H$ has vertex set $\ZZ\times V(H)$, with
\[
  (i,a)\sim(i,b)\ \text{ for } ab\in E(H),\qquad (i,a)\sim(i+1,a)\ \text{ for } a\in V(H).
\]
With the columns $V_i=\{i\}\times V(H)$ and the shift $\sigma(i,a)=(i+1,a)$,
the graph $\widetilde G$ is periodic with span $h=1$.
For $H=K_2$ it is the infinite ladder of Example~\ref{ex:ladder}.
For $a\in V(H)$ we call $\ZZ\times\{a\}$ a \emph{track}, which induces a two-way infinite path.

\begin{figure}[ht]
  \centering
  \begin{tikzpicture}[x=9mm, y=7mm,
      vertex/.style = {circle, fill=white, draw=black, line width=.5pt,
                       inner sep=0pt, minimum size=1.8mm},
      edge/.style   = {line width=.5pt, black!60}]
    % top: the cover P_infty x P_3
    \begin{scope}
      \foreach \y in {0,1,2} {
        \draw[edge] (-0.55,\y) -- (6.55,\y);
        \node at (-0.95,\y) {$\cdots$};
        \node at (6.95,\y) {$\cdots$};
      }
      \foreach \i in {0,...,6} {
        \draw[edge] (\i,0) -- (\i,2);
        \foreach \y in {0,1,2} \node[vertex] at (\i,\y) {};
        \pgfmathtruncatemacro{\c}{\i-1}
        \node[font=\small] at (\i,-0.85) {$\c$};
      }
      \node[vertex, fill=black] at (1,2) {};
      \node[vertex, fill=black] at (6,0) {};
      \node[font=\small] at (1,2.5) {$x$};
      \node[font=\small, anchor=north west] at (6.02,-0.05) {$\tau(x)$};
      \node[font=\small, anchor=west] at (7.4,1) {$P_\infty\,\square\,P_3$};
    \end{scope}
    % bottom: the twisted cylinder C_5 x_phi P_3
    \begin{scope}[yshift=-3.1cm, xshift=9mm]
      \foreach \y in {0,1,2} \draw[edge] (-0.5,\y) -- (4.5,\y);
      \foreach \i in {0,...,4} {
        \draw[edge] (\i,0) -- (\i,2);
        \foreach \y in {0,1,2} \node[vertex] at (\i,\y) {};
        \node[font=\small] at (\i,-0.8) {$\i$};
      }
      \foreach \y/\l in {2/\alpha, 1/\beta, 0/\gamma} \node[font=\small] at (4.85,\y) {$\l$};
      \node[vertex, fill=black] at (0,2) {};
      \node[font=\small] at (0,2.5) {$\pi(x)$};
      \foreach \y/\l in {2/\gamma, 1/\beta, 0/\alpha} \node[font=\small] at (-0.85,\y) {$\l$};
      \node[font=\small, anchor=west] at (6.4,1) {$C_5\,\square_\varphi P_3$};
    \end{scope}
  \end{tikzpicture}
  \caption{The cover $P_\infty\,\square\,P_3$ (top)
    and the twisted cylinder $C_5\,\square_\varphi P_3$ with $\varphi$ the reflection of $P_3$ (bottom),
    where half-edges with the same label are joined.}
  \label{fig:cylinder}
\end{figure}
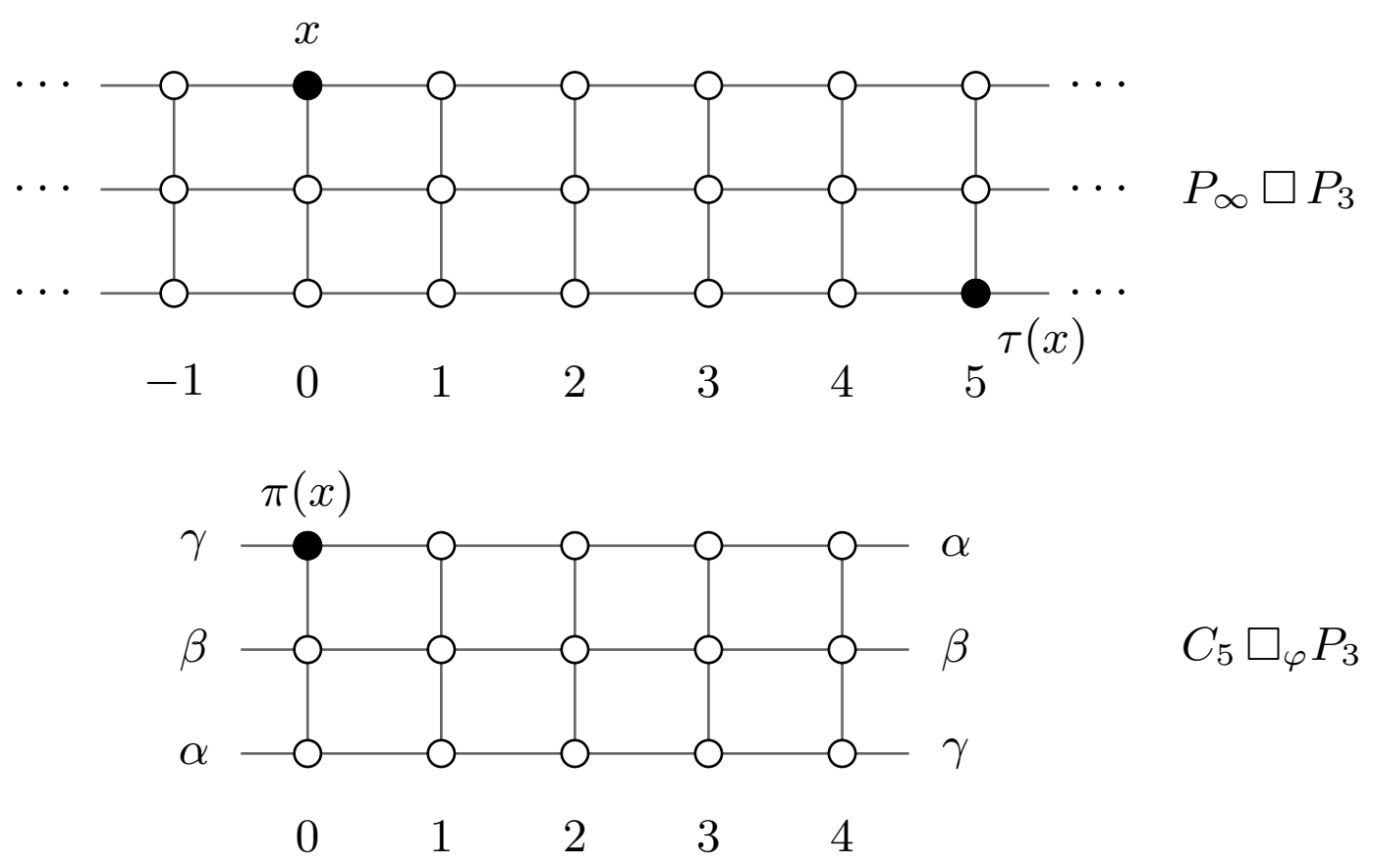

Let $\varphi$ be an automorphism of $H$ and $n\ge 3$.
The \emph{twisted cylinder} $C_n\,\square_\varphi H$ has vertex set $\{0,\ldots,n-1\}\times V(H)$,
in which $(i,a)$ is adjacent to $(i,b)$ for $ab\in E(H)$,
$(i,a)$ is adjacent to $(i+1,a)$ for $0\le i<n-1$,
and $(n-1,a)$ is adjacent to $(0,\varphi(a))$.
For $\varphi=1$ this is the Cartesian product $C_n\,\square\, H$,
and for $H=K_2$ with $\varphi$ the nontrivial automorphism it is the Möbius ladder on $2n$ vertices, as in Example~\ref{ex:ladder}.
The map $\tau(i,a)=(i+n,\varphi^{-1}(a))$ is an automorphism of $\widetilde G$ and a shift by $n$.
Writing $i=qn+r$ with $0\le r<n$,
the map $\pi(i,a)=(r,\varphi^{q}(a))$ is constant on the orbits of $\langle\tau\rangle$
and identifies $\widetilde G/\langle\tau\rangle$ with $C_n\,\square_\varphi H$.
See Figure~\ref{fig:cylinder}.

\begin{lemma}\label{lem:two_layers}
  Two consecutive columns form a zero forcing set of $\widetilde G$.
\end{lemma}

\begin{proof}
  Suppose $V_i\cup V_{i+1}$ is black.
  Each vertex $(i+1,a)$ has the black neighbors $(i,a)$ and $(i+1,b)$ with $ab\in E(H)$,
  so it forces its only other neighbor $(i+2,a)$, and $V_{i+2}$ becomes black.
  Likewise each $(i,a)$ forces $(i-1,a)$, and $V_{i-1}$ becomes black.
  Repeating this with $V_{i+1}\cup V_{i+2}$ and $V_{i-1}\cup V_i$, every column becomes black.
\end{proof}

\begin{lemma}\label{lem:spread}
  The graph $\widetilde G$ has the $1$-window property.
\end{lemma}

\begin{proof}
  We apply Lemma~\ref{lem:window_criterion} to $W=V_i$.
  It is connected since $H$ is connected,
  and $N[V_i]=V_{i-1}\cup V_i\cup V_{i+1}$ is a zero forcing set by Lemma~\ref{lem:two_layers}.
\end{proof}

\begin{proposition}\label{prop:Z_cyl_inf}
  Let $H$ be a finite connected graph with $m$ vertices. Then
  \[
    \Z(P_\infty\,\square\, H)=2m.
  \]
\end{proposition}

\begin{proof}
  The upper bound follows from Lemma~\ref{lem:two_layers}.
  The $m$ tracks are disjoint two-way infinite paths, so the lower bound follows from Corollary~\ref{cor:tracks}.
\end{proof}

\begin{theorem}\label{thm:cylinder}
  Let $H$ be a finite connected graph with $m$ vertices, and let $\varphi$ be an automorphism of $H$.
  For all $n\ge\max\{2m,3\}$,
  \[
    \Z(C_n\,\square_\varphi H)=2m.
  \]
  In particular, $\Z(C_n\,\square\, H)=2m$ for all $n\ge\max\{2m,3\}$.
\end{theorem}

\begin{proof}
  We apply Theorem~\ref{thm:periodic_Z} with $h=g=1$.
  Condition~(a) holds since $\widetilde G$ is connected and each column has $m$ vertices,
  (b) holds by Lemma~\ref{lem:spread},
  and (c) holds since $\Z(\widetilde G)=2m$ by Proposition~\ref{prop:Z_cyl_inf}.
  The threshold in Theorem~\ref{thm:periodic_Z} is
  $g(\Z(\widetilde G)-2)+h+1=(2m-2)+1+1=2m$.
  Since $n\ge 3>2h$ and $n\ge 2m$,
  Theorem~\ref{thm:periodic_Z} gives $\Z(C_n\,\square_\varphi H)=2m$.
\end{proof}

By~\cite[Corollary~2.8]{AIM2008}, $\Z(C_n\,\square\, H)\le\min\{\Z(H)\,n,\,2m\}$.
For $n\ge 2m$ this bound equals $2m$ since $\Z(H)\ge 1$,
and Theorem~\ref{thm:cylinder} shows that it is attained.
The threshold cannot be lowered in general:
for $H=P_m$ the bound is $\min\{n,2m\}$, which is the exact value by~\cite[Theorem~3.8]{AIM2008}
and is less than $2m$ for $n<2m$.
Several cases of Theorem~\ref{thm:cylinder} are known, usually with a smaller threshold for the particular $H$:
the cylinders $C_n\,\square\, P_m$, including the prisms $C_n\,\square\, K_2$~\cite[Theorem~3.8]{AIM2008},
the Möbius ladders~\cite[Proposition~3.9]{AIM2008},
the tori $C_n\,\square\, C_m$~\cite[Theorem~2.6]{Benson2018},
and $C_n\,\square\, K_m$ and its twist by a cyclic permutation~\cite[Theorem~4.1]{Duong2020}.
Most of these results also show $\M=\Z$, which is stronger since $\M\le\Z$ for every graph.
These proofs construct matrices of the required nullity or use the structure of the particular $H$.
General lower bounds for Cartesian products, such as $\Z(G\,\square\, G')\ge\Z(G)\Z(G')+1$
under a condition on fractional zero forcing~\cite{CHKMS2023}, apply to arbitrary factors and all $n$.
On the products $C_n\,\square\, H$ with $n\ge\max\{2m,3\}$ they give at most $2\Z(H)+1\le 2m-1$,
while Theorem~\ref{thm:cylinder} gives the exact value $2m$.

\section{Circulants}\label{sec:circulants}

\begin{figure}[ht]
  \centering
  \begin{tikzpicture}[x=8mm, y=8mm,
      vertex/.style = {circle, fill=white, draw=black, line width=.5pt,
                       inner sep=0pt, minimum size=1.8mm},
      edge/.style   = {line width=.5pt, black!60}]
    % top: the cover C_infty(1,3)
    % vertex i is drawn in column i and row i mod 3
    \begin{scope}
      \begin{scope}
        \clip (-0.6,-0.5) rectangle (11.6,2.5);
        \foreach \y in {0,1,2} \draw[edge] (-1,\y) -- (12,\y);
        \foreach \i in {-1,...,11} {
          \pgfmathtruncatemacro{\r}{mod(\i+3,3)}
          \pgfmathtruncatemacro{\t}{mod(\i+4,3)}
          \draw[edge] (\i,\r) -- (\i+1,\t);
        }
      \end{scope}
      \foreach \y in {0,1,2} {
        \node at (-1,\y) {$\cdots$};
        \node at (12,\y) {$\cdots$};
      }
      \foreach \i in {0,...,11} {
        \pgfmathtruncatemacro{\r}{mod(\i,3)}
        \node[vertex] at (\i,\r) {};
        \node[font=\small] at (\i,-0.7) {$\i$};
      }
      \node[vertex, fill=black] at (1,1) {};
      \node[vertex, fill=black] at (11,2) {};
      \node[font=\small, anchor=south east] at (0.95,1.05) {$x$};
      \node[font=\small] at (11,2.5) {$\tau(x)$};
      \node[font=\small, anchor=west] at (12.5,1) {$C_\infty(1,3)$};
    \end{scope}
    % bottom: the circulant C_10(1,3)
    \begin{scope}[yshift=-3.2cm, xshift=4.4cm]
      \foreach \i in {0,...,9} {
        \coordinate (c\i) at ({90-36*\i}:1.6cm);
      }
      \foreach \i in {0,...,9} {
        \pgfmathtruncatemacro{\j}{mod(\i+1,10)}
        \pgfmathtruncatemacro{\l}{mod(\i+3,10)}
        \draw[edge] (c\i) -- (c\j);
        \draw[edge] (c\i) -- (c\l);
      }
      \foreach \i in {0,...,9} {
        \node[vertex] at (c\i) {};
        \node[font=\small] at ({90-36*\i}:1.95cm) {$\i$};
      }
      \node[vertex, fill=black] at (c1) {};
      \node[font=\small] at (54:2.5cm) {$\pi(x)$};
      \node[font=\small, anchor=west] at (2.6,0) {$C_{10}(1,3)$};
    \end{scope}
  \end{tikzpicture}
  \caption{The cover $C_\infty(1,3)$ (top) and the circulant $C_{10}(1,3)$ (bottom).}
  \label{fig:circulant}
\end{figure}
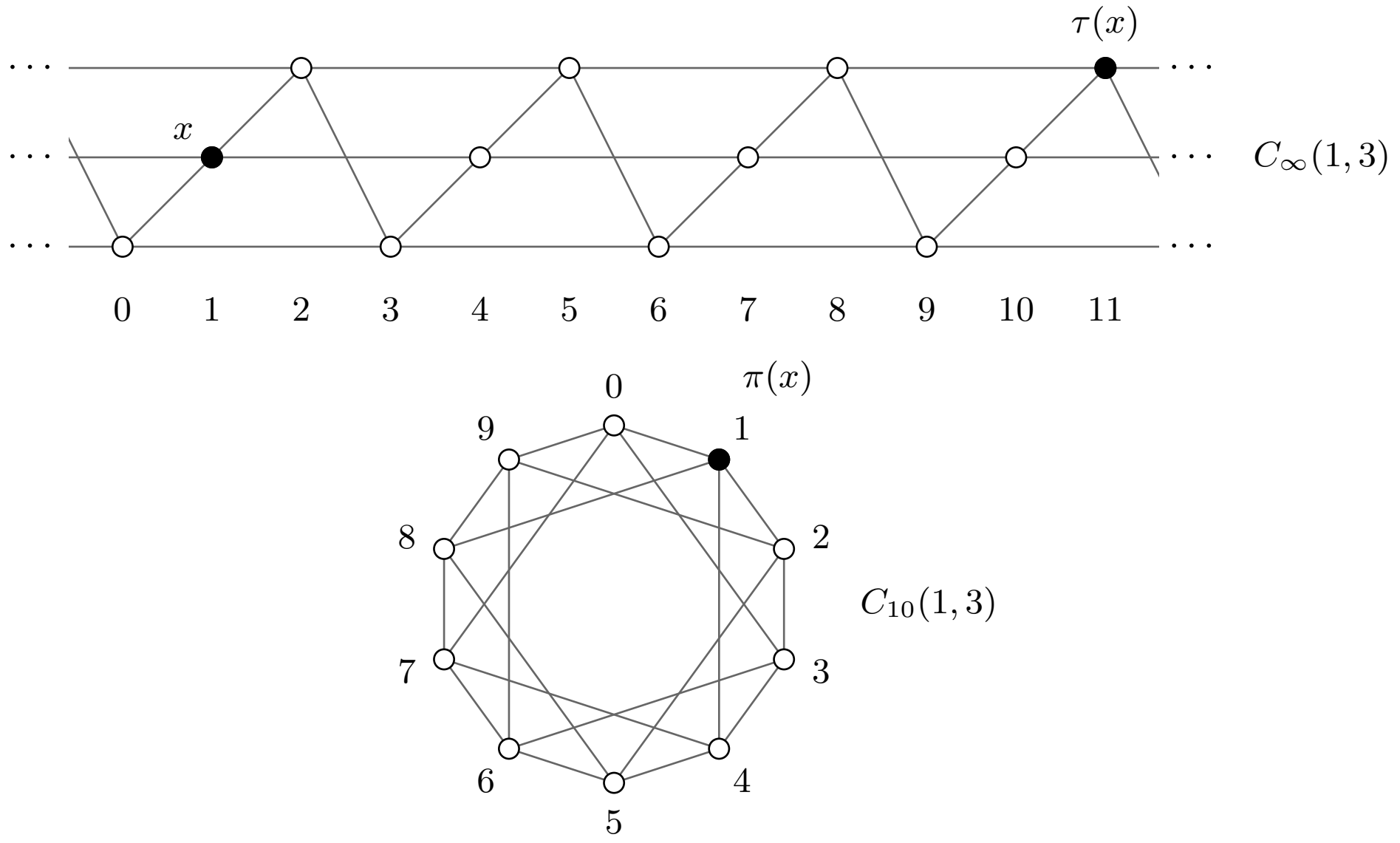

Fix an integer $a\ge 2$.
For $n>2a$, the circulant $C_n(1,a)$ has vertex set $\ZZ/n\ZZ$, with $i\sim i\pm 1$ and $i\sim i\pm a$.
Its cover $\widetilde G=C_\infty(1,a)$ has vertex set $\ZZ$ and the same adjacencies.
Since $a\ge 2$, the graph $\widetilde G$ is $4$-regular, with
\[
  N(i)=\{i-a,\,i-1,\,i+1,\,i+a\}.
\]
With the columns $\{i\}$ and the translation $\sigma(i)=i+1$,
the graph $\widetilde G$ is periodic with span $h=a$, and each column is a single vertex.
The map $\tau(i)=i+n$ is an automorphism of $\widetilde G$ and a shift by $n$,
and $\pi(i)=i\bmod n$ identifies $\widetilde G/\langle\tau\rangle$ with $C_n(1,a)$.
See Figure~\ref{fig:circulant}.
The vertices of a residue class modulo $a$ induce a two-way infinite path in $\widetilde G$.
These $a$ paths play the role of the tracks.

\begin{lemma}\label{lem:circ_consecutive}
  Any $2a$ consecutive vertices form a zero forcing set of $\widetilde G$.
\end{lemma}

\begin{proof}
  By translation, we may assume $1,2,\ldots,2a$ are black.
  The vertex $a+1$ has the black neighbors $a$, $a+2$ and $1$,
  so it forces its fourth neighbor $2a+1$.
  Likewise $a$ has the black neighbors $a-1$, $a+1$ and $2a$, so it forces $0$.
  Repeating this on both sides, every vertex becomes black.
\end{proof}

\begin{lemma}\label{lem:circ_window}
  The graph $\widetilde G$ has the $a$-window property.
\end{lemma}

\begin{proof}
  We apply Lemma~\ref{lem:window_criterion} to $W=\{0,\ldots,a-1\}$.
  It is connected by the edges $x\sim x+1$,
  and $N[W]=\{-a,\ldots,2a-1\}$ contains $2a$ consecutive vertices,
  so it is a zero forcing set by Lemma~\ref{lem:circ_consecutive}.
\end{proof}

\begin{proposition}\label{prop:circ_inf}
  $\Z(C_\infty(1,a))=2a$.
\end{proposition}

\begin{proof}
  The upper bound is Lemma~\ref{lem:circ_consecutive}.
  The $a$ residue classes modulo $a$ induce disjoint two-way infinite paths,
  so the lower bound is Corollary~\ref{cor:tracks}.
\end{proof}

\begin{theorem}\label{thm:circulant}
  For $a\ge 2$ and $n\ge 2a^2-a+1$, $\Z(C_n(1,a))=2a$.
\end{theorem}

\begin{proof}
  We apply Theorem~\ref{thm:periodic_Z} with $h=g=a$.
  Condition~(a) holds since $\widetilde G$ is connected and each column is a single vertex,
  (b) holds by Lemma~\ref{lem:circ_window},
  and (c) holds since $\Z(\widetilde G)=2a$ by Proposition~\ref{prop:circ_inf}.
  The threshold in Theorem~\ref{thm:periodic_Z} is
  $g(\Z(\widetilde G)-2)+h+1=a(2a-2)+a+1=2a^2-a+1$.
  Since $n\ge 2a^2-a+1>2a$,
  Theorem~\ref{thm:periodic_Z} gives $\Z(C_n(1,a))=2a$.
\end{proof}

As for $\GP{n}{k}$, there is a least integer $N'_a$ such that $\Z(C_n(1,a))=2a$ for all $n\ge N'_a$.

\begin{proposition}\label{prop:N_a}
  For every $a\ge 3$,
  \[
    a^2+1\le N'_a\le 2a^2-a+1,
  \]
  and $N'_2=5$.
\end{proposition}

\begin{proof}
  The upper bound is Theorem~\ref{thm:circulant}.
  For the lower bound, $\Z(C_{a^2}(1,a))\le 2a-1$ when $a\ge 3$ by~\cite[Corollary~4.11]{Duong2020},
  so $N'_a>a^2$.
  For $a=2$, $C_n(1,2)$ is a consecutive circulant, so $\Z(C_n(1,2))=4$ for all $n\ge 5$,
  the least $n$ allowed, by~\cite[Theorem~5.4]{DeaettMeyer2016}; see also Table~\ref{tab:circ}.
\end{proof}

We determined $\Z(C_n(1,a))$ for $2a<n<2a^2-a+1$ and $a\le 6$ by an exhaustive computer search.
The results are in Tables~\ref{tab:circ} and~\ref{tab:circ6}.
They show $N'_a=a^2+1$ for $2\le a\le 6$, so the lower bound of Proposition~\ref{prop:N_a} is attained in these cases.
The values agree with $\Z(C_9(1,3))=5$ and $\Z(C_{3m}(1,3))=6$ for $m>3$
from~\cite[Theorems~4.12 and~4.15]{Duong2020}.
At $n=2a+1$ the value is $4$, since multiplication by $2$ is an isomorphism $C_{2a+1}(1,a)\cong C_{2a+1}(1,2)$.

\begin{table}[ht]
  \caption{The zero forcing number $\Z(C_n(1,a))$ for $2\le a\le 5$.}
  \label{tab:circ}
  \setlength{\tabcolsep}{2.4pt}
  \begin{tabular}{c|*{21}{c}|c|c}
    $a\backslash n$ & 5 & 6 & 7 & 8 & 9 & 10 & 11 & 12 & 13 & 14 & 15 & 16 & 17 & 18 & 19 & 20 & 21 & 22 & 23 & 24 & 25 & $\ge 26$ & $N'_a$ \\
    \hline
    2 & 4 & 4 & 4 & 4 & 4 & 4 & 4 & 4 & 4 & 4 & 4 & 4 & 4 & 4 & 4 & 4 & 4 & 4 & 4 & 4 & 4 & 4 & 5 \\
    3 &   &   & 4 & 6 & 5 & 6 & 6 & 6 & 6 & 6 & 6 & 6 & 6 & 6 & 6 & 6 & 6 & 6 & 6 & 6 & 6 & 6 & 10 \\
    4 &   &   &   &   & 4 & 6 & 6 & 6 & 6 & 6 & 7 & 7 & 8 & 8 & 8 & 8 & 8 & 8 & 8 & 8 & 8 & 8 & 17 \\
    5 &   &   &   &   &   &   & 4 & 8 & 6 & 6 & 6 & 6 & 6 & 8 & 8 & 8 & 8 & 8 & 8 & 10 & 9 & 10 & 26 \\
  \end{tabular}
\end{table}

\begin{table}[ht]
  \caption{The zero forcing number $\Z(C_n(1,6))$, which is constant on each range of $n$ shown.}
  \label{tab:circ6}
  \centering
  \begin{tabular}{c|*{7}{c}|c}
    $a\backslash n$ & 13 & 14 & 15--20 & 21--27 & 28--34 & 35--36 & $\ge 37$ & $N'_a$ \\
    \hline
    6 & 4 & 8 & 6 & 8 & 10 & 11 & 12 & 37 \\
  \end{tabular}
\end{table}

We close with a question about the thresholds.
By Propositions~\ref{prop:N_k} and~\ref{prop:N_a}, both $N_k$ and $N'_a$ are quadratic,
but the leading coefficients of the upper and lower bounds differ by a factor of $2$.
The known values $N_3=13$, $N_4=18$, $N_5=23$ and $N'_a=a^2+1$ for $2\le a\le 6$
are much closer to the lower bounds than to the upper bounds.

\noindent\textbf{Question.}
Can the upper bounds $2k^2+k+1$ and $2a^2-a+1$ be lowered to $k^2+O(k)$ and $a^2+O(a)$?

\bibliographystyle{amsplain}
\bibliography{paper}

\end{document}